\documentclass[a4, 11pt, leqno]{article}
\usepackage{amsmath}
\usepackage{amssymb}
\newtheorem{df}{Definition}[section]
\newtheorem{pr}[df]{Proposition}
\newtheorem{Th}[df]{Theorem}

\newtheorem{lm}[df]{Lemma}

\begin{document}
\title{\bf Motion of a Vortex Filament in the Half Space}
\author{Masashi A{\sc iki} and Tatsuo I{\sc guchi}}
\date{}
\maketitle
\vspace*{-0.5cm}
\begin{center}
Department of Mathematics, Faculty of Science and Technology, Keio University, \\
3-14-1 Hiyoshi, Kohoku-ku, Yokohama, 223-8522, J{\sc apan}
\end{center}

\begin{abstract}
A model equation for the motion of a vortex filament immersed in three dimensional, 
incompressible and inviscid fluid is investigated as a humble attempt to model the 
motion of a tornado. 
We solve an initial-boundary value problem in the half space where we impose a boundary 
condition in which the vortex filament is allowed to move on the boundary. 
\end{abstract}

\section{Introduction}
Many researchers have studied tornadoes from several perspectives. 
A systematic research and observation of tornadoes is difficult mainly because of two reasons: 
a precise prediction of tornado formation is not yet possible, 
and the life-span of a tornado is very short, giving only short openings for any kind of measurements. 
Many aspects of tornadoes are still unknown. 

In 1971, Fujita \cite{1} gave a systematic categorization of tornadoes. 
He proposed the so-called Fujita scale in which tornadoes are classified according to 
the damage that it dealt to buildings and other surroundings. 
The scale provides a correlation between ranges of wind speed and the damage that it causes. 
The enhanced version, called the Enhanced Fujita Scale, is used to classify tornadoes to date. 
McDonald \cite{8} gives a review of Fujita's contributions to tornado research. 

Since then, due to the advancement of technology, more accurate and thorough observations 
and simulations have become possible, and theories for the formation and motion of tornadoes have developed. 
Klemp \cite{2} and the references within give an extensive review on the known dynamics of tornadoes. 

Motivated by this, we investigate the motion of a vortex filament. 
The vortex filament equation, 
also called the Localized Induction Equation (LIE) models the movement of a vortex filament, 
which is a space curve where the vorticity of the fluid is concentrated, and is described by 
\begin{eqnarray}
\mbox{\mathversion{bold}$x$}_{t} = \mbox{\mathversion{bold}$x$}_{s}\times \mbox{\mathversion{bold}$x$}_{ss}, 
\label{intro}
\end{eqnarray}
where
\( \mbox{\mathversion{bold}$x$}(s,t)=\left( x^{1}(s,t) , \ x^{2}(s,t), \ x^{3}(s,t) \right) \) 
is the position of the vortex filament parameterized by the arc length \( s\) at time \( t\), 
\( \times \) denotes the exterior product, and the subscripts denote differentiation with respect to that variable. 
We also use \( \partial _{s} \) and \( \partial _{t} \) for partial differentiation with the corresponding variables. 

The LIE was first derived by Da Rios \cite{3} and re-examined by Arms and Hama \cite{9}. 
Since then, many authors have worked with the equation. 
Nishiyama and Tani \cite{10, 4} gave the unique solvability of the initial value problem for the LIE 
in Sobolev spaces. 
A different approach was taken by Hasimoto \cite{6}. 
He used the so-called Hasimoto transformation to transform (\ref{intro}) into a nonlinear Schr\"odinger equation: 
\begin{equation*}
\frac{1}{\rm i}\frac{\partial \psi}{\partial t}
= \frac{ \partial ^{2}\psi }{\partial s^{2}} + \frac{1}{2} \left| \psi \right| ^{2}\psi ,
\end{equation*}
where \( \psi \) is given by 
\[ \psi = \kappa \exp \left( {\rm i} \int ^{s}_{0} \tau \,{\rm d}s \right), \]
\( \kappa \) is the curvature, and \( \tau \) is the torsion of the filament. 
Even though this expression is undefined at points of the filament where the curvature vanishes, 
Koiso \cite{7} proved that the Hasimoto transformation is well-defined in the class of \( C^{\infty}\) 
functions. 
He used a geometrical approach to define the Hasimoto transformation and showed the unique solvability 
of the initial value problem in the class of \( C^{\infty}\) functions. 

Regarding initial-boundary value problems, 
the only known result that the authors know is by Nishiyama and Tani \cite{4}. 
The boundary condition imposed there necessarily fixes the end point of the vortex filament 
and does not allow it to move on the boundary. 
From the physical point of view, the vortex filament must be closed, extend to the spatial infinity, 
or end on boundaries of the fluid region. 
In the last case, we have to impose an appropriate boundary condition to show the well-posedness of the problem. 
Since it is hard to find what kind of boundary condition is physically reasonable if we begin our analysis from 
the Schr\"odinger equation, we chose to work with the original vortex filament equation (\ref{intro}). 

In light of modeling the motion of a tornado, we consider (\ref{intro}) in a framework 
in which the end of the vortex filament is allowed to move on the boundary. 
We do this by setting a different boundary condition than that of \cite{4}. 

The contents of this paper are as follows. 
In section 2, we formulate our problem and give basic notations. 
In section 3, we derive compatibility conditions for our initial-boundary value problem. 
Sections 4 and 5 are concerned with constructing the solution.

\section{Setting of the Problem}
\setcounter{equation}{0}

We consider the initial-boundary value problem for the motion of a vortex filament in the half-space 
in which the filament is allowed to move on the boundary: 
\begin{equation}
\left\{
 \begin{array}{lll}
  \mbox{\mathversion{bold}$x$}_{t}=\mbox{\mathversion{bold}$x$}_{s}
   \times \mbox{\mathversion{bold}$x$}_{ss}, & & s>0, \ t>0, \\
  \mbox{\mathversion{bold}$x$}(s,0)=\mbox{\mathversion{bold}$x$}_{0}(s),  & & s>0, \\
  \mbox{\mathversion{bold}$x$}_{s}(0,t)=\mbox{\mathversion{bold}$e$}_{3}, & & t>0,
 \end{array}
\right.
\label{p1}
\end{equation}
where \( \mbox{\mathversion{bold}$e$}_{3}= (0,0,1) \). 
We assume that 
\begin{equation}
|\mbox{\mathversion{bold}$x$}_{0s}(s)|=1 \makebox[3em]{for} s\geq 0, \qquad
x^{3}_{0}(0)=0,
\label{a1}
\end{equation}
for the initial datum. 
The first condition states that the initial vortex filament is parametrized by the arc length and 
the second condition just states that the curve is parameterized starting from the boundary. 
Here we observe that by taking the inner product of \( \mbox{\mathversion{bold}$e$}_{3} \) with the equation, 
taking the trace at \( s=0\), and noting the boundary condition we have 
\begin{equation*}
\begin{aligned}
\frac{d}{dt}\left( \mbox{\mathversion{bold}$e$}_{3}\cdot \mbox{\mathversion{bold}$x$}\right) |_{s=0} 
 &= \left. \mbox{\mathversion{bold}$e$}_{3}\cdot \left( \mbox{\mathversion{bold}$x$}_{s}\times 
\mbox{\mathversion{bold}$x$}_{ss}\right) \right|_{s=0} \\
 &= \left. \mbox{\mathversion{bold}$x$}_{s}\cdot \left( \mbox{\mathversion{bold}$x$}_{s}\times 
  \mbox{\mathversion{bold}$x$}_{ss}\right) \right|_{s=0} \\
 &=0,
\end{aligned}
\end{equation*}
where \( ``\cdot" \) denotes the inner product and \( |_{s=0} \) denotes the trace at \( s=0 \). 
This means that if the end of the vortex filament is on the boundary initially, 
then it will stay on the boundary, but is not necessarily fixed. 
This is our reason for the notion ``allowed to move on the boundary".

By introducing new variables
\( \mbox{\mathversion{bold}$v$}(s,t) := \mbox{\mathversion{bold}$x$}_{s}(s,t) \) and 
\( \mbox{\mathversion{bold}$v$}_{0}(s) := \mbox{\mathversion{bold}$x$}_{0s}(s) \), 
(\ref{p1}) and (\ref{a1}) become 
\begin{equation}
\left\{
 \begin{array}{lll}
  \mbox{\mathversion{bold}$v$}_{t}=\mbox{\mathversion{bold}$v$}
   \times \mbox{\mathversion{bold}$v$}_{ss}, & & s>0, \ t>0, \\
  \mbox{\mathversion{bold}$v$}(s,0)=\mbox{\mathversion{bold}$v$}_{0}(s), & & s>0, \\
  \mbox{\mathversion{bold}$v$}(0,t)=\mbox{\mathversion{bold}$e$}_{3}, & & t>0, 
 \end{array}
\right.
\label{p2}
\end{equation}
\begin{equation}
|\mbox{\mathversion{bold}$v$}_{0}(s)|=1, \quad s\geq 0.
\label{a2}
\end{equation}
Once we solve (\ref{p2}), the solution \( \mbox{\mathversion{bold}$x$} \) of (\ref{p1}) and (\ref{a1}) 
can be constructed by 
\[ \mbox{\mathversion{bold}$x$}(s,t)
 = \mbox{\mathversion{bold}$x$}_{0}(s)+ \int ^{t}_{0} \mbox{\mathversion{bold}$v$}(s,\tau ) \times 
 \mbox{\mathversion{bold}$v$}_{s}(s,\tau ) \,{\rm d}\tau. \]
So from now on, we concentrate on the initial-boundary value problem (\ref{p2}) 
under the condition (\ref{a2}). 
Note that if the initial datum satisfies (\ref{a2}), then any smooth solution 
\( \mbox{\mathversion{bold}$v$} \) of (\ref{p2}) satisfies 
\begin{equation}
|\mbox{\mathversion{bold}$v$}(s,t)|=1, \quad s\geq 0 , \ t\geq 0.
\end{equation}
This can be confirmed by taking the inner product of the equation with \( \mbox{\mathversion{bold}$v$} \).

\bigskip
We define basic notations that we will use throughout this paper. 

For a domain \( \Omega \), 
a non-negative integer \( m\), and \( 1\leq p \leq \infty \), \( W^{m,p}(\Omega )\) is the Sobolev space 
containing all real-valued functions that have derivatives in the sense of distribution up to order \( m\) 
belonging to \( L^{p}(\Omega )\). 
We set \( H^{m}(\Omega ) = W^{m,2}(\Omega ) \) as the Sobolev space equipped with the usual inner product, 
and \( H^{0}(\Omega ) = L^{2}(\Omega ) \). 
We will particularly use the cases \( \Omega= {\mathbf R} \) and \( \Omega = {\mathbf R}_{+} \), 
where \( {\mathbf R}_{+} = \{ s\in {\mathbf R}; s>0 \} \). 
The norm in \( H^{m}(\Omega ) \) is denoted by \( ||\cdot ||_{m} \) and we simply write \( ||\cdot || \) 
for \( ||\cdot ||_{0} \). 
We do not indicate the domain in the symbol for the norms 
since we use it in a way where there is no risk of confusion. 

For a Banach space \( X\), 
\( C^{m}([0,T];X) \) denotes the spaces of functions that are \( m\) times continuously differentiable 
in \( t\) with respect to the topology of \( X\).

For any function space described above, we say that a vector valued function belongs to the function space 
if each of its components does.

\section{Compatibility Conditions}
\setcounter{equation}{0}
We derive necessary conditions for a smooth solution to exist for (\ref{p2}) with (\ref{a2}). 

Suppose that \( \mbox{\mathversion{bold}$v$}(s,t) \) is a smooth solution of (\ref{p2}) with (\ref{a2}) 
defined in \( {\mathbf R}_{+}\times [ 0,T ] \) for some positive \( T \). 
We have already seen that for all \( (s,t) \in {\mathbf R}_{+}\times [0,T] \) 
\begin{equation}
|\mbox{\mathversion{bold}$v$}(s,t)|^{2} =1.
\label{u}
\end{equation}
By differentiating the boundary condition with respect to \( t \) we see that
\[  \left. \partial _{t}^{n}\mbox{\mathversion{bold}$v$}\right| _{s=0} = \mbox{\mathversion{bold}$0$} 
 \makebox[3.5em]{\rm for} n\in {\mathbf N}, \ t>0. 
 \leqno{(B)_{n}} \]
We next show 

\begin{lm}
For a smooth solution \( \mbox{\boldmath $v$}(s,t)\) under consideration, 
it holds that
\[ \left. \mbox{\mathversion{bold}$v$}\times \partial_{s}^{2n}\mbox{\mathversion{bold}$v$} \right| _{s=0}
 = \mbox{\mathversion{bold}$0$}, 
 \leqno{(C)_{n}} \]
\vspace{-6ex}
\[ \left. \partial_{s}^{j}\mbox{\mathversion{bold}$v$} \cdot \partial _{s}^{l}\mbox{\mathversion{bold}$v$} 
 \right| _{s=0} =0 \makebox[3.5em]{\rm for} j+l=2n+1. 
 \leqno{(D)_{n}} \]
\label{lemma}
\end{lm}

\medskip
\noindent
{\it Proof}. 
We prove them by induction. 
From \( (B)_{1} \) and by taking the trace of the equation we see that 
\[ \mbox{\mathversion{bold}$0$} = \mbox{\mathversion{bold}$v$}_{t}\left. \right| _{s=0}
 = \mbox{\mathversion{bold}$v$}\times \mbox{\mathversion{bold}$v$}_{ss}\left. \right| _{s=0},\]
thus, \( (C)_{1}\) holds. 
By taking the exterior product of \( \mbox{\mathversion{bold}$v$}_{s}\) and (\( C\))\(_{1}\) we have 
\[ \left\{ \left( \mbox{\mathversion{bold}$v$}_{s}\cdot \mbox{\mathversion{bold}$v$}_{ss} \right) 
 \mbox{\mathversion{bold}$v$} - 
\left( \mbox{\mathversion{bold}$v$}_{s}\cdot \mbox{\mathversion{bold}$v$}\right) 
 \mbox{\mathversion{bold}$v$}_{ss}\right\} \left. \right| _{s=0} = \mbox{\mathversion{bold}$0$}. \] 
On the other hand, by differentiating (\ref{u}) with respect to \( s\) we have 
\( \mbox{\mathversion{bold}$v$}\cdot \mbox{\mathversion{bold}$v$}_{s} \equiv 0 \). 
Combining these two and the fact that \( \mbox{\mathversion{bold}$v$}\) is a non-zero vector, 
we arrive at 
\[ \mbox{\mathversion{bold}$v$}_{s}\cdot \mbox{\mathversion{bold}$v$}_{ss}\left. \right| _{s=0} =0. \]
Finally, by differentiating (\ref{u}) with respect to \( s\) three times and setting \( s=0 \), we have
\[ 0= 2\left( \mbox{\mathversion{bold}$v$}\cdot \mbox{\mathversion{bold}$v$}_{sss}
 + 3\mbox{\mathversion{bold}$v$}_{s}\cdot \mbox{\mathversion{bold}$v$}_{ss} \right) \left. \right|_{s=0}
  = 2\mbox{\mathversion{bold}$v$}\cdot \mbox{\mathversion{bold}$v$}_{sss}\left. \right| _{s=0}, \]
so, \( (D)_{1} \) holds. 

Suppose that the statements hold up to \( n-1 \) for some \( n\geq 2 \). 
By differentiating \( (C)_{n-1} \) with respect to \( t\) we have 
\[ \mbox{\mathversion{bold}$v$}\times \bigl( 
 \partial ^{2(n-1)}_{s}\mbox{\mathversion{bold}$v$}_{t} \bigr)  \bigr| _{s=0} 
 = \mbox{\mathversion{bold}$0$}, \]
where we have used \( (B)_{1} \). 
We see that 
\begin{eqnarray*}
\partial^{2(n-1)}_{s} \mbox{\mathversion{bold}$v$}_{t}
 = \partial^{2(n-1)}_s \left( \mbox{\mathversion{bold}$v$}\times \mbox{\mathversion{bold}$v$}_{ss}\right)
 = \sum ^{2(n-1)}_{k=0} 
\left(
 \begin{array}{c}
  2(n-1)\\
  k
 \end{array}
\right)
\bigl( \partial ^{k}_{s}\mbox{\mathversion{bold}$v$}\times
 \partial ^{2(n-1)-k+2}_{s}\mbox{\mathversion{bold}$v$} \bigr),
\end{eqnarray*}
where 
\(
\left(
 \begin{array}{c}
  2(n-1)\\
  k
 \end{array}
\right)
\)
is the binomial coefficient. 
So we have 
\begin{eqnarray}
\sum ^{2(n-1)}_{k=0}
\left(
 \begin{array}{c}
  2(n-1)\\
  k
 \end{array}
\right)
\left\{ \mbox{\mathversion{bold}$v$}\times 
 \bigl( \partial ^{k}_{s}\mbox{\mathversion{bold}$v$}\times
 \partial ^{2(n-1)-k+2}_{s}\mbox{\mathversion{bold}$v$} \bigr)\right\}\biggr| _{s=0} 
 = \mbox{\mathversion{bold}$0$}. \label{1}
\end{eqnarray}
We examine each term in the summation. 
When \( 2\leq k \leq 2(n-1) \) is even, we see from the assumptions of induction \( (C)_{k/2} \) and 
\( (C)_{(2(n-1)-k+2)/2} \) that both \( \partial ^{k}_{s}\mbox{\mathversion{bold}$v$} \) and 
\( \partial ^{2(n-1)-k+2}_{s}\mbox{\mathversion{bold}$v$} \) are parallel to 
\( \mbox{\mathversion{bold}$v$} \), so that 
\[ \left. \partial ^{k}_{s}\mbox{\mathversion{bold}$v$}\times
 \partial ^{2(n-1)-k+2}_{s}\mbox{\mathversion{bold}$v$}\right| _{s=0} = \mbox{\mathversion{bold}$0$}. \]
When \( 1\leq k \leq 2(n-1) \) is odd, we rewrite the exterior product in (\ref{1}) as 
\begin{eqnarray*}
\mbox{\mathversion{bold}$v$}\times \bigl( \partial ^{k}_{s}\mbox{\mathversion{bold}$v$}\times
 \partial ^{2(n-1)-k+2}_{s}\mbox{\mathversion{bold}$v$} \bigr)
= \bigl( \mbox{\mathversion{bold}$v$}\cdot \partial ^{2(n-1)-k+2}_{s}\mbox{\mathversion{bold}$v$} \bigr)
 \partial ^{k}_{s}\mbox{\mathversion{bold}$v$}
 - \bigl( \mbox{\mathversion{bold}$v$}\cdot \partial ^{k}_{s}\mbox{\mathversion{bold}$v$} \bigr) 
 \partial ^{2(n-1)-k+2}_{s}\mbox{\mathversion{bold}$v$}.
\end{eqnarray*}
Since \( 2(n-1)-k +2 \) is also odd, by \( (D)_{(k-1)/2} \ {\rm and} \ (D)_{(2(n-1)-k+1)/2} \) we have 
\[ \left. \mbox{\mathversion{bold}$v$}\cdot \partial ^{k}_{s}\mbox{\mathversion{bold}$v$} \right| _{s=0}
 = \left. \mbox{\mathversion{bold}$v$}\cdot \partial ^{2(n-1)-k+2}_{s}\mbox{\mathversion{bold}$v$} \right| _{s=0}
 =0. \]
Thus, only the term with \( k=0 \) remains and we get 
\[ \left. \mbox{\mathversion{bold}$v$}\times \bigl( \mbox{\mathversion{bold}$v$}\times
 \partial ^{2n}_{s}\mbox{\mathversion{bold}$v$} \bigr) \right| _{s=0} = \mbox{\mathversion{bold}$0$}. \]
Here, we note that 
\[ \mbox{\mathversion{bold}$v$}\times (\mbox{\mathversion{bold}$v$}\times
 \partial ^{2n}_{s}\mbox{\mathversion{bold}$v$} ) 
= (\mbox{\mathversion{bold}$v$}\cdot \partial ^{2n}_{s}\mbox{\mathversion{bold}$v$} )\mbox{\mathversion{bold}$v$}
 - \partial ^{2n}_{s}\mbox{\mathversion{bold}$v$}, \]
where we used (\ref{u}). 
Taking the exterior product of this with \( \mbox{\mathversion{bold}$v$}\) we see that \( (C)_{n} \) holds. 
Taking the exterior product of \( \partial ^{2n+1-2k}_{s}\mbox{\mathversion{bold}$v$} \) with \( (C)_{k} \) 
and using \( (D)_{n-k} \) for \( 1\leq k \leq n \) yields 
\[ \left. \bigl( \partial ^{2k}_{s}\mbox{\mathversion{bold}$v$} \cdot
 \partial ^{2n+1-2k}_{s}\mbox{\mathversion{bold}$v$} \bigr) \mbox{\mathversion{bold}$v$} \right| _{s=0} 
 = \mbox{\mathversion{bold}$0$}. \]
Since \( \mbox{\mathversion{bold}$v$}\) is non zero, we have for \( 1\leq k \leq n \) 
\begin{eqnarray}
\left. \partial ^{2k}_{s}\mbox{\mathversion{bold}$v$} \cdot
 \partial^{2n+1-2k}_{s}\mbox{\mathversion{bold}$v$} \right| _{s=0} = 0. \label{2}
\end{eqnarray}
Finally, by differentiating (\ref{u}) with respect to \( s\) \( (2n+1) \) times, we have 
\[ \sum ^{2n+1}_{j=0} 
\left(
 \begin{array}{c}
  2n+1\\
  j
 \end{array}
\right)
\left( \partial ^{j}_{s}\mbox{\mathversion{bold}$v$} \cdot
 \partial ^{2n+1-j}_{s}\mbox{\mathversion{bold}$v$} \right) \Bigr| _{s=0} =0. \]
Since every term except \( j=0, \ 2n+1 \) is of the form (\ref{2}), we see that 
\[ \left. \mbox{\mathversion{bold}$v$}\cdot \partial ^{2n+1}_{s}\mbox{\mathversion{bold}$v$} \right| _{s=0} =0, \]
which, together with (\ref{2}), finishes the proof of \( (D)_{n} \). 
\hfill \( \Box \)

\bigskip
Worth noting are the following two properties which will be used in later parts of this paper. 
For integers \( n\), 
\[ \left. \mbox{\mathversion{bold}$e$}_{3}\times \partial ^{2n}_{s}\mbox{\mathversion{bold}$v$}\right| _{s=0}
 = \mbox{\mathversion{bold}$0$}, 
 \quad \left. \mbox{\mathversion{bold}$e$}_{3}\cdot \partial ^{2n+1}_{s}\mbox{\mathversion{bold}$v$}\right| _{s=0}=0.
\]
These are special cases of \( (C)_{n}\) and \( (D)_{n}\) with the boundary condition substituted in. 

By taking the limit \( t\rightarrow 0 \) in \( (C)_{n} \), we derive a necessary condition for the initial datum.

\begin{df}
For \( n\in {\mathbf N}\cup \{ 0\} \), we say that the initial datum \( v_{0} \) satisfies 
the compatibility condition \( (A)_{n}\) if the following condition is satisfied for \( 0\leq k \leq n \) 
\[
\left\{
 \begin{array}{lll}
  \left. \mbox{\mathversion{bold}$v$}_{0}\right| _{s=0}=\mbox{\mathversion{bold}$e$}_{3}, & & k=0, \\[1ex]
  \left. \left( \mbox{\mathversion{bold}$v$}_{0}\times
   \partial ^{2k}_s \mbox{\mathversion{bold}$v$}_{0}\right) \right|_{s=0}=\mbox{\mathversion{bold}$0$}, 
   & & k\in {\mathbf N}.
 \end{array}
\right.
\]
\end{df}

\medskip
From the proof of Lemma \ref{lemma}, we see that if \( \mbox{\mathversion{bold}$v$}_{0} \) satisfies (\ref{a2}) 
and the compatibility condition \( (A)_{n} \), then \( \mbox{\mathversion{bold}$v$}_{0} \) also satisfies 
\( (D)_{k} \) for \( 0\leq k\leq n \) with \( \mbox{\mathversion{bold}$v$} \) replaced by 
\( \mbox{\mathversion{bold}$v$}_{0} \) as long as the trace exists.

\section{Extension of the Initial Datum}
\setcounter{equation}{0}
For the initial datum \( \mbox{\mathversion{bold}$v$}_{0} \) defined on the half-line, 
we extend it to the whole line by 
\begin{equation}
\widetilde{\mbox{\mathversion{bold}$v$}}_{0}(s) = 
\left\{
 \begin{array}{lll}
  \mbox{\mathversion{bold}$v$}_{0}(s), & & s\geq 0, \\[.5ex]
  -\overline{\mbox{\mathversion{bold}$v$}}_{0}(-s), & & s<0,
 \end{array}
\right.
\label{extention}
\end{equation}
where \( \overline{\mbox{\mathversion{bold}$v$}} = (v^{1},v^{2},-v^{3}) \) for 
\( \mbox{\mathversion{bold}$v$}= (v^{1},v^{2},v^{3}) \in {\mathbf R}^{3}\).

\begin{pr}
For any integer \ \( m\geq 2 \), if \( \mbox{\mathversion{bold}$v$}_{0s} \in H^{m}({\mathbf R}_{+}) \) 
satisfies {\rm (\ref{a2})} and the compatibility condition \( (A)_{[ \frac{m}{2} ]} \), then 
\( \widetilde{\mbox{\mathversion{bold}$v$}}_{0s} \in H^{m}({\mathbf R}) \). 
Here, \( [ \frac{m}{2} ] \) indicates the largest integer not exceeding \( \frac{m}{2} \). 
\label{sp2}
\end{pr}

\noindent
{\it Proof}. 
Fix an arbitrary integer \( m \geq 2 \). We will prove by induction on \( k\) that 
\( \partial _{s}^{k+1}\widetilde{\mbox{\mathversion{bold}$v$}}_{0} \in L^{2}({\mathbf R}) \) 
for any \( 0\leq k \leq m \). 
Specifically we show that the derivatives of \( \widetilde{\mbox{\mathversion{bold}$v$}}_{0}\) 
in the distribution sense on the whole line \( {\mathbf R}\) up to order \( m+1\) have the form 
\begin{equation}
\bigl( \partial ^{k+1}_{s}\widetilde{\mbox{\mathversion{bold}$v$}}_{0}\bigr) (s) = 
\left\{
 \begin{array}{lll}
  \left( \partial ^{k+1}_{s}\mbox{\mathversion{bold}$v$}_{0}\right) (s), & & s>0, \\[1ex]
  -(-1)^{k+1} \bigl( \overline{ \partial ^{k+1}_{s}\mbox{\mathversion{bold}$v$}_{0}} \bigr) (-s), & & s<0,
 \end{array}
\right.
\label{derivative}
\end{equation}
for \( 0\leq k\leq m\). 

Since \( \mbox{\mathversion{bold}$v$}_{0}\in L^{\infty}({\mathbf R}_{+})\) and 
\( \mbox{\mathversion{bold}$v$}_{0s}\in H^{2}({\mathbf R}_{+}) \), 
Sobolev's embedding theorem states \( \mbox{\mathversion{bold}$v$}_{0s}\in L^{\infty}({\mathbf R}_{+}) \) 
and thus \( \mbox{\mathversion{bold}$v$}_{0}\in W^{1,\infty}({\mathbf R}_{+}) \), 
so that the trace \( \mbox{\mathversion{bold}$v$}_{0}(0) \) exists. 
By definition (\ref{extention}) we have 
\[ \widetilde{\mbox{\mathversion{bold}$v$}}_{0}(-0) = \left( -v_{0}^{1}(0), \ -v_{0}^{2}(0), \ v^{3}(0) \right),\]
but from \( (A)_{0} \), \( v_{0}^{1}(0) = v_{0}^{2}(0) =0 \). 
These imply that \( \widetilde{\mbox{\mathversion{bold}$v$}}_{0}(+0)
=\widetilde{\mbox{\mathversion{bold}$v$}}_{0}(-0) \), so that we obtain 
\begin{eqnarray*}
\partial _{s}\widetilde{\mbox{\mathversion{bold}$v$}}_{0}(s) = 
\left\{
 \begin{array}{lll}
  \left( \partial _{s}\mbox{\mathversion{bold}$v$}_{0}\right) (s), & & s>0, \\[1ex]
  - (-1)\bigl( \overline{ \partial _{s}\mbox{\mathversion{bold}$v$}_{0}} \bigr) (-s), & & s<0 ,
 \end{array}
\right.
\end{eqnarray*}
and the case \( k=0 \) is proved. 

Suppose that (\ref{derivative}) with \(k+1\) replaced by \(k\) holds for some \( k\in \{ 1,2,\ldots ,m\} \). 
We check that the derivative \( \partial^{k}_{s}\widetilde{\mbox{\mathversion{bold}$v$}}_{0} \) 
does not have a jump discontinuity at \( s=0\). 
When \( k\) is even, from the definition of \( \overline{\partial ^{k}_{s}\mbox{\mathversion{bold}$v$}_{0}} \), 
\[ \bigl( \partial^{k}_{s}\widetilde{\mbox{\mathversion{bold}$v$}}_{0}\bigr) (-0)
 = \bigl( -\partial ^{k}_{s}v_{0}^{1}(0), 
  -\partial ^{k}_{s}v_{0}^{2}(0), \ \partial ^{k}_{s}v_{0}^{3}(0) \bigr), \]
but from \( (A)_{\frac{k}{2}} \) we have 
\[ \mbox{\mathversion{bold}$0$} = \left. \mbox{\mathversion{bold}$v$}_{0}\times
 \partial ^{k}_{s}\mbox{\mathversion{bold}$v$}_{0}\right| _{s=0} = 
\mbox{\mathversion{bold}$e$}_{3}\times \partial ^{k}_{s}\mbox{\mathversion{bold}$v$}_{0}(0), \]
which means that \( \partial ^{k}_{s}\mbox{\mathversion{bold}$v$}_{0}(0) \) is parallel to 
\( \mbox{\mathversion{bold}$e$}_{3} \) and that the first and second components are zero. 
When \( k \) is odd, 
\[ \bigl( \partial^{k}_{s}\widetilde{\mbox{\mathversion{bold}$v$}}_{0}\bigr) (-0)
 =  \bigl( \partial ^{k}_{s}v_{0}^{1}(0),\ 
\partial ^{k}_{s}v_{0}^{2}(0), \ -\partial ^{k}_{s}v_{0}^{3}(0) \bigr), \]
but \( (A)_{[\frac{k}{2}]} \) implies \( (D)_{[\frac{k}{2}]} \) and particularly 
\[ 0= \left. \mbox{\mathversion{bold}$v$}_{0}\cdot
 \partial ^{k}_{s}\mbox{\mathversion{bold}$v$}_{0} \right| _{s=0}
 =\mbox{\mathversion{bold}$e$}_{3}\cdot \partial ^{k}_{s}\mbox{\mathversion{bold}$v$}_{0}(0) 
 = (\partial _{s}^{k}v_{0}^{3})(0), \]
so the third component is zero. 
In both cases, we have \( \bigl(\partial ^{k}_{s}\widetilde{\mbox{\mathversion{bold}$v$}}_{0}\bigr) (+0)
= \bigl(\partial ^{k}_{s}\widetilde{\mbox{\mathversion{bold}$v$}}_{0}\bigr) (-0) \), so that 
we can verify (\ref{derivative}). 
This finishes the proof of the proposition. 
\hfill \( \Box \)

\section{Existence and Uniqueness of Solution}
\setcounter{equation}{0}
Using \( \widetilde{\mbox{\mathversion{bold}$v$}}_{0} \), we consider the following initial value problem: 
\begin{eqnarray}
&& \mbox{\mathversion{bold}$u$}_{t} = \mbox{\mathversion{bold}$u$}\times \mbox{\mathversion{bold}$u$}_{ss}, 
 \qquad s\in {\mathbf R} , \ t>0,  \label{p3}\\
&& \mbox{\mathversion{bold}$u$}(s,0) = \widetilde{\mbox{\mathversion{bold}$v$}}_{0}(s),
 \quad\: s\in {\mathbf R}. \label{p4}
\end{eqnarray}
By Proposition \ref{sp2}, the existence and uniqueness theorem (cf. Nishiyama \cite{5}) of a strong solution 
\( \mbox{\mathversion{bold}$u$}\) is applicable. 
Specifically we use the following theorem.

\begin{Th}{\rm (Nishiyama \cite{5})}
For a non-negative integer m, if \( \widetilde{\mbox{\mathversion{bold}$v$}}_{0s}\in H^{2+m}({\mathbf R}) \) and 
\( |\widetilde{\mbox{\mathversion{bold}$v$}}_{0}| \equiv 1 \), then the initial value problem {\rm (\ref{p3})} and 
{\rm (\ref{p4})} has a unique solution \mbox{\mathversion{bold}$u$} such that 
\[ \mbox{\mathversion{bold}$u$}-\widetilde{\mbox{\mathversion{bold}$v$}}_{0}
 \in C\bigl( [0,\infty ); H^{3+m}({\mathbf R})\bigr) \cap C^{1}\bigl( [0,\infty ); H^{1+m}({\mathbf R})\bigr) \]
and \( |\mbox{\mathversion{bold}$u$}|\equiv 1 \).
\end{Th}

From Proposition \ref{sp2}, the assumptions of the theorem are satisfied if 
\( \mbox{\mathversion{bold}$v$}_{0s}\in H^{2+m}({\mathbf R}_{+}) \) satisfies the compatibility condition 
\( (A)_{[\frac{2+m}{2}]} \) and (\ref{a2}). 

Now we define the operator T by 
\[ ({\rm T}\mbox{\mathversion{bold}$w$})(s) = -\overline{\mbox{\mathversion{bold}$w$}}(-s), \]
for \( {\mathbf R}^{3} \)-valued functions \( \mbox{\mathversion{bold}$w$}\) defined on \( s\in {\mathbf R} \). 
By direct calculation, we can verify that 
\( {\rm T}\widetilde{\mbox{\mathversion{bold}$v$}}_{0} = \widetilde{\mbox{\mathversion{bold}$v$}}_{0} \) 
and that 
\( {\rm T}\left( \mbox{\mathversion{bold}$u$}\times \mbox{\mathversion{bold}$u$}_{ss} \right)
 = \left( {\rm T}\mbox{\mathversion{bold}$u$}\right)
\times \left( {\rm T}\mbox{\mathversion{bold}$u$}\right) _{ss}. \)
Taking these into account and applying the operator T to (\ref{p3}) and (\ref{p4}), we have 
\begin{eqnarray*}
\left\{
 \begin{array}{lll}
  ({\rm T}\mbox{\mathversion{bold}$u$})_{t}
   = ({\rm T}\mbox{\mathversion{bold}$u$}) \times ({\rm T}\mbox{\mathversion{bold}$u$}) _{ss},
   & & s\in {\mathbf R} , \ t>0, \\[.5ex]
  ({\rm T}\mbox{\mathversion{bold}$u$})(s,0) = ({\rm T}\widetilde{\mbox{\mathversion{bold}$v$}}_{0})(s)
   = \widetilde{\mbox{\mathversion{bold}$v$}}_{0}(s),
   & & s\in {\mathbf R}, 
\end{array}\right.
\end{eqnarray*}
which means that T\( \mbox{\mathversion{bold}$u$}\) is also a solution of (\ref{p3}) and (\ref{p4}). 
Thus we have T\( \mbox{\mathversion{bold}$u$}\) = \( \mbox{\mathversion{bold}$u$}\) 
by the uniqueness of the solution. 
Therefore, for any \( t\in [0,T] \) 
\[ \mbox{\mathversion{bold}$u$}(0,t) = \left( {\rm T}\mbox{\mathversion{bold}$u$}\right) (0,t)
 = -\overline{\mbox{\mathversion{bold}$u$}}(0,t), \]
which is equivalent to \( u^{1}(0,t) = u^{2}(0,t) =0 \). 
Therefore, it holds that \( u^{3}(0,t)= -1 \ {\rm or} \ 1 \) because we have 
\( |\mbox{\mathversion{bold}$u$}| \equiv 1 \). 
But in view of \( \widetilde{\mbox{\mathversion{bold}$v$}}_{0}(0) = \mbox{\mathversion{bold}$v$}_{0}(0)
 = \mbox{\mathversion{bold}$e$}_{3} \), we obtain 
\( \mbox{\mathversion{bold}$u$}(0,t) = \mbox{\mathversion{bold}$e$}_{3} \) by the continuity in \( t \).

This shows that the restriction of \( \mbox{\mathversion{bold}$u$}\) to \( {\mathbf R}_{+} \) 
is a solution of our initial-boundary value problem. 
Using this function \( \mbox{\mathversion{bold}$v$} := \mbox{\mathversion{bold}$u$}|_{\mathbf{R}_{+}}\), 
we can construct the solution \( \mbox{\mathversion{bold}$x$}\) to the original equation 
as we stated in section 2. 
Thus we have

\begin{Th}
For a non-negative integer \( m\), if \( \mbox{\mathversion{bold}$x$}_{0ss} \in H^{2+m}({\mathbf R}_{+}) \) 
and \( x_{0s} \) satisfies the compatibility condition \( (A)_{[\frac{2+m}{2}]} \) and {\rm (\ref{a1})}, 
then there exists a unique solution \( \mbox{\mathversion{bold}$x$}\) of {\rm (\ref{p1})} such that 
\[ \mbox{\mathversion{bold}$x$}-\mbox{\mathversion{bold}$x$}_{0}
 \in C\bigl( [0,\infty ); H^{4+m}({\mathbf R}_{+})\bigr) \cap
  C^{1}\bigl( [0,\infty ); H^{2+m}({\mathbf R}_{+})\bigr), \]
and \( |\mbox{\mathversion{bold}$x$}_{s}|\equiv 1 \). 
\end{Th}

\noindent
{\it Proof.} 
The uniqueness is left to be proved. 
Suppose that \( \mbox{\mathversion{bold}$x$}_{1}\) and \( \mbox{\mathversion{bold}$x$}_{2}\) 
are solutions as in the theorem. 
Then, by extending \( \mbox{\mathversion{bold}$x$}_{i} \ (i=1,2)\) by 
\begin{eqnarray*}
\widetilde{\mbox{\mathversion{bold}$x$}}_{i}(s,t)=
\left\{
 \begin{array}{lll}
  \mbox{\mathversion{bold}$x$}_{i}(s,t) & & s\geq 0, t>0, \\
  \overline{\mbox{\mathversion{bold}$x$}}_{i}(-s,t) & & s<0,t>0,
 \end{array}
\right.
\end{eqnarray*}
we see that \( \widetilde{\mbox{\mathversion{bold}$x$}}_{i} \) are solutions of the vortex filament equation 
in the whole space. 
Thus \( \mbox{\mathversion{bold}$x$}_{1}=\mbox{\mathversion{bold}$x$}_{2}\) follows from the uniqueness 
of the solution in the whole space.
\hfill \( \Box \)


\end{document}